\documentclass[12pt]{amsart}

\copyrightinfo{2019}{}

\usepackage{latexsym, amsmath, amssymb,amsthm,amscd,amsopn,amsfonts,mathrsfs,tikz}
\usepackage{stmaryrd}
\usepackage{version}
\usepackage{epsfig,graphics,color,graphicx,graphpap,dsfont,eufrak}
\usepackage{amssymb}
\usepackage{tikz}
\usetikzlibrary{snakes}

\ifx\pdfoutput\undefined
  \DeclareGraphicsExtensions{.eps}
\else
  \ifx\pdfoutput\relax
    \DeclareGraphicsExtensions{.eps}
  \else
    \ifnum\pdfoutput>0
      \DeclareGraphicsExtensions{.pdf}
    \else
      \DeclareGraphicsExtensions{.eps}
    \fi
  \fi
\fi

\usepackage{subfigure}

\newtheorem{theorem}{Theorem}[section]

\theoremstyle{definition}

\theoremstyle{remark}

\numberwithin{equation}{section}

\newcommand{\ga}{\gamma}

\newcommand{\dl}{\delta}
\newcommand{\Dl}{\Delta}

\newcommand{\ra}{\rightarrow}

\newcommand{\sg}{\sigma}

\newcommand{\pa}{\partial}

\newcommand{\La}{\Lambda}

\newcommand{\om}{\omega}
\newcommand{\Om}{\Omega}
\newcommand{\na}{\nabla}

\newcommand{\non}{\nonumber}

\begin{document}

\title[Infinite norm of the derivative of the solution operator]{Infinite norm of the derivative of the solution operator of Euler equations}

\author{Y. Charles Li}
\address{Department of Mathematics, University of Missouri, 
Columbia, MO 65211, USA}
\email{liyan@missouri.edu}
\urladdr{http://faculty.missouri.edu/~liyan}

\curraddr{}
\thanks{}

\subjclass{Primary 76, 35; Secondary 34}

\date{}

\dedicatory{}


\begin{abstract}
Through a simple and elegant argument, we prove that the norm of the derivative of the solution operator of Euler equations posed in the Sobolev space $H^n$, along any base solution that is in $H^n$ but not in $H^{n+1}$, is infinite. We also review the counterpart of this result for 
Navier-Stokes equations at high Reynolds number from the perspective of fully developed turbulence. Finally we present 
a few examples and numerical simulations to show a more complete picture of the so-called rough dependence upon initial data.
\end{abstract}

\maketitle

\tableofcontents

\section{Introduction}

The solution operator of Euler equations of fluids is nowhere differentiable \cite{Inc15} \cite{IL18}. This is what we called ``rough dependence on initial data" for Euler equations.
There are several ways for the solution operator to be non-differentiable. The most common way is that the norm of the derivative of the solution operator is infinite. Other ways include that the norm of the formal derivative of the solution operator is finite but the Fr\'echet definition of derivative is violated. Our Main Theorem of this paper is:
\begin{theorem}
If the Euler equations are posed in the Sobolev space $H^n$ and the base solution is in $H^n$ but not in $H^{n+1}$,
then the norm of the formal derivative of the solution operator is infinite when $t>0$.
\label{MT}
\end{theorem}

If $u_0$ is an element in $H^n$ but not in $H^{n+1}$, then for any $v_0$ in $H^{n+1}$, $u_0 + v_0$ is an element in $H^n$ but not in $H^{n+1}$. Thus there are more elements in $H^n$ but not in $H^{n+1}$ than in $H^{n+1}$. When $n > \frac{d}{2} +1$ ($d$ is the spatial dimension), each initial element in $H^n$ generates a local solution of the Euler equations in $H^n$, and each initial element in $H^{n+1}$ generates a local solution of the Euler equations in $H^{n+1}$. Thus Theorem \ref{MT} is valid for a majority of base solutions in $H^n$. 

Even though it is everywhere differentiable, the solution operator of Navier-Stokes equations will somehow approach the
solution operator of Euler equations in the infinite Reynolds number limit. This is the regime that
we are most interested in from the perspective of fully developed turbulence. In this regime,
we believe that the norm of the derivative of the solution operator along turbulent solutions of Navier-Stokes equations will approach infinity in the infinite Reynolds number limit. Since the norm of the derivative of the solution operator measures the maximal growth of perturbations, perturbations in
fully developed turbulence grow superfast. This is what we called ``rough dependence on initial data" for fully developed turbulence. Our theory is that fully developed turbulence is initiated and maintained by such superfast growth of ever existing perturbations \cite{Li14}. Such superfast
growth can reach substantial scale even in short time, and leads to superfast nonlinear saturation and short term unpredictability
of fully developed turbulence. The superfast growth of perturbations also implies that the turbulent solution of Navier-Stokes equations and the turbulent flow in reality (lab or 
nature) can be substantially different in short time, even though they have the same initial condition. 

\section{Basic formulation}

The  Navier-Stokes equations are given by
\begin{eqnarray}
& & u_t - \frac{1}{Re} \Dl u  = - \na p - u\cdot \na u , \label{NS1} \\
& & \na \cdot u = 0 ,  \label{NS2}
\end{eqnarray}
where $u$ is the $d$-dimensional fluid velocity ($d=2,3$), $p$ is the fluid pressure, and 
$Re$ is the Reynolds number. Setting the Reynolds number to infinity $Re = \infty$, the 
Navier-Stokes equations (\ref{NS1})-(\ref{NS2}) reduce to the Euler equations
\begin{eqnarray}
& & u_t  = - \na p - u\cdot \na u , \label{E1} \\
& & \na \cdot u = 0 . \label{E2}
\end{eqnarray}
For any $u \in H^n(\mathbb{R}^d)$  ($n > \frac{d}{2} +1$), there 
is a neighborhood $B$ and a short time $T>0$, such that for any $v \in B$ there exists a unique 
solution to the Navier-Stokes equations (\ref{NS1})-(\ref{NS2}) in $C^0([0,T]; H^n(\mathbb{R}^d))$. As $Re 
\ra \infty$, this solution converges to that of the Euler equations (\ref{E1})-(\ref{E2}) in the same space \cite{Kat72} \cite{Kat75}. For any 
$t \in [0, T]$, let $S^t$ be the solution map:
\begin{equation}
S^t \  :  \   B  \mapsto H^n(\mathbb{R}^d), \  S^t (u(0)) = u(t), \label{SM} 
\end{equation}
i.e. the solution map maps the initial condition to the solution's value at time $t$. The solution map 
is continuous for both Navier-Stokes equations  (\ref{NS1})-(\ref{NS2}) and Euler equations (\ref{E1})-(\ref{E2}) \cite{Kat72} 
\cite{Kat75}, but nowhere differentiable for  Euler equations \cite{Inc15} \cite{IL18}. Even though the derivative of the solution map for Navier-Stokes equations (\ref{NS1})-(\ref{NS2}) exists, it is natural to conjecture that the norm of the derivative of the solution 
map along turbulent solutions approaches infinity as the Reynolds number approaches infinity.  The following upper bound was 
obtained in \cite{Li14}.
\begin{equation}
\| DS^t(u(0)) \| = \sup_{d u(0)} \frac{\| d u(t) \|_n}{\| d u(0) \|_n} \leq e^{\sg \sqrt{Re} \sqrt{t} \ + \ \sg_1 t}, \label{UB} 
\end{equation}
where $d u(0)$ is any initial perturbation of $u(0)$, $\| \  \|_n$ represents the Sobolev $H^n$ norm, and 
\[
\sg = \frac{8c}{\sqrt{2e}} \max_{\tau \in [0,T]} \| u(\tau )\|_n, \ \ \sg_1 = \frac{\sqrt{2e}}{2} \sg ,
\]
where $c$ only depends on $n$ and the spatial domain.
The above bound also applies to spatially periodic domain $\mathbb{T}^d$ instead of $\mathbb{R}^d$.

Sometimes, it is convenient to use the Leray projection of the Navier-Stokes equations. The Leray projection is an orthogonal projection in $L^2(\mathbb{R}^d)$, given by
\[
\mathbb{P} g = g - \na \Dl^{-1} \na \cdot g . 
\]
Applying the Leray projection to the Navier-Stokes equations  (\ref{NS1})-(\ref{NS2}), one gets
\begin{equation}
u_t + \frac{1}{Re} \Dl u  = - \mathbb{P} \left ( u\cdot \na u \right ) , \label{NSL} 
\end{equation}
and the corresponding Euler equations 
\begin{equation}
u_t = - \mathbb{P} \left ( u\cdot \na u \right ) . \label{EL} 
\end{equation}

\section{Proof of the Main Theorem}

Here we will present a simple and elegant proof of the Main Theorem \ref{MT}. We will present the periodic boundary condition case, of course the same proof applies to other boundary condition cases that allow translational invariance.

\begin{proof}
The Euler equations (\ref{E1})-(\ref{E2}) are translationally invariant, i.e. if $u(t,x)$ is a solution, then $u(t,x-at)+a$ are also solutions for constant vectors $a$. Using 
Fourier series, we have
\begin{eqnarray*}
&& u(t,x) = \sum_{k\in \mathbb{Z}^d} u_k(t) e^{ik\cdot x}, \\
&& U(t,x,a) = u(t,x-at)+a = a + \sum_{k\in \mathbb{Z}^d} u_k(t) e^{ik\cdot (x-at)} .
\end{eqnarray*}
Let $u(t,x)$ be any solution that is in $H^n$ but not in $H^{n+1}$, then $U(t,x,a)$ is a family of solutions parametrized by $a$, which have the same property, and 
$U(t,x,0)=u(t,x)$. Notice also that 
\[
U(0,x,a) = u(0,x)+a.
\]
By varying $a$ around $a=0$, we can make a directional variation of the initial condition around $u(0,x)$, which leads to a directional derivative of the solution operator:
\[
\frac{\pa}{\pa a_m}U(t,x,a) |_{a=0} = \frac{\pa}{\pa a_m} a + \sum_{k\in \mathbb{Z}^d} (-ik_mt)u_k(t) e^{ik\cdot x}, \  m=1, \cdots, d.
\]
Thus
\begin{eqnarray*}
&& \sum_{m=1}^d \| \frac{\pa}{\pa a_m}U(t,x,a) |_{a=0}  \|_n^2 \\
&& = d (2\pi )^d + t^2 (2\pi )^d \sum_{k\in \mathbb{Z}^d} \left ( |k|^2 + \cdots + |k|^{2(n+1)}\right ) |u_k|^2 \\
&& = d (2\pi )^d + t^2 \left ( \| u(t,x)\|^2_{n+1} - \| u(t,x)\|^2_0 \right ).
\end{eqnarray*}
Since $\| u(t,x)\|_{n+1} = \infty$, 
\[
\| \frac{\pa}{\pa a_m}U(t,x,a) |_{a=0}  \|_n = \infty , \quad \text{when } t>0 ,
\]
for some $m$. Since it is the supremum over the norms of all directional derivatives, the norm of the derivative of the solution operator along $u(t,x)$ is infinite, and 
this completes the proof of the main theorem.
\end{proof}

\section{Example 1 - the trivial one}

Under either periodic condition or decaying boundary condition in the whole space, the base
solution is the trivial solution $u=0$ ($p=0$). Then the corresponding linearized
equations of (\ref{NSL}) are given by
\[
du_t - \frac{1}{Re} \Dl du = 0,
\]
and the corresponding linearized equations of (\ref{EL}) are given by
\[
du_t = 0.
\]
Thus
\begin{equation}
du(t) = e^{\frac{t}{Re} \Dl} du(0).
\label{Ex11}
\end{equation}
Starting from the same initial condition $\dl u(0) = du(0)$, the increment $\dl u(t)$ satisfies
\[
\dl u_t - \frac{1}{Re} \Dl \dl u = -\mathbb{P} (\dl u \cdot \na \dl u),
\]
in the Navier-Stokes case, and
\[
\dl u_t =  -\mathbb{P} (\dl u \cdot \na \dl u),
\]
in the Euler case. By the method of variation of parameters,
\begin{equation}
\dl u(t) = e^{\frac{t}{Re} \Dl} \dl u(0) - \int_0^t e^{\frac{t-\tau }{Re} \Dl}
\mathbb{P} (\dl u \cdot \na \dl u) d\tau ,
\label{Ex12}
\end{equation}
in the Navier-Stokes case. Applying the inequality
\[
\| e^{\frac{t}{Re} \Dl} u \|_n \leq \left ( \frac{1}{\sqrt{2e}} \sqrt{\frac{Re}{t}} +1 \right )
\| u \|_{n-1},
\]
where $\| \ \|_n$ represents the Sobolev $H^n$ norm ($n > \frac{d}{2} + 1$),
\[
\|\dl u(t)\|_n \leq \| \dl u(0)\|_n + 2 c\int_0^t  \left ( \frac{1}{\sqrt{2e}}
\sqrt{\frac{Re}{t - \tau }} +1 \right ) \|\dl u(\tau )\|^2_n d\tau ,
\]
where $c$ is a constant that only depends on $n$ and the spatial domain. Then
\[
\max_{t\in [0, T]}\|\dl u(t)\|_n \leq \| \dl u(0)\|_n + 2 \max_{t\in [0, T]}\|\dl u(t)\|^2_n 
\int_0^t  \left ( \frac{1}{\sqrt{2e}} \sqrt{\frac{Re}{t - \tau }} +1 \right ) d\tau .
\]
Thus
\[
\max_{t\in [0, T]}\|\dl u(t)\|_n \sim \| \dl u(0)\|_n, \quad \text{as } \| \dl u(0)\|_n \ra 0.
\]
In view of the fact that $\dl u(0) = du(0)$, (\ref{Ex12})-(\ref{Ex11}) leads to
\[
\dl u(t) - du(t) = - \int_0^t e^{\frac{t-\tau }{Re} \Dl}
\mathbb{P} (\dl u \cdot \na \dl u) d\tau ,
\]
which can be estimated as above,
\[
\| \dl u(t) - du(t) \|_n \sim \| \dl u(0)\|^2_n, \quad \text{as } \| \dl u(0)\|_n \ra 0.
\]
Thus when $Re < \infty$, the derivative of the solution operator at the trivial solution exists, and
is given by (\ref{Ex11}). But in the Euler case,
\[
\dl u(t) - du(t) = - \int_0^t \mathbb{P} (\dl u \cdot \na \dl u) d\tau ,
\]
and $\| \dl u(t) - du(t) \|_n$ is not of order $o(\| \dl u(0)\|_n)$, as $\| \dl u(0)\|_n \ra 0$.
Thus in the Euler case, the derivative of the solution operator at the trivial solution still
does not exist, even though the norm of the formal derivative is bounded. 

\section{Example 2 - the simple one}

The base solution is the 2D Couette linear shear $u=(x_2, 0)$. The boundary conditions are
\[
u = (\pm 1, 0), \quad \text{at } x_2 = \pm 1,
\]
in the viscous case, no $u_1$ condition (slip) in the inviscid case, and periodic boundary
condition along $x_1$-direction with period $2\pi$. It is more convenient to use the vorticity
variable
\[
\om = \frac{\pa u_2}{\pa x_1} -  \frac{\pa u_1}{\pa x_2} = - \Dl \psi ,
\]
where
\[
u_1 = \frac{\pa \psi}{\pa x_2}, \ u_2 = -\frac{\pa \psi}{\pa x_1}.
\]
In terms of the vorticity variable, the 2D Navier-Stokes equations take the form
\[
\om_t - \frac{1}{Re} \Dl \om = - u \cdot \na \om .
\]
The linearized 2D Navier-Stokes equation at the Couette linear shear is given by
\[
d\om_t - \frac{1}{Re} \Dl d\om = - x_2 d\om_{x_1} .
\]
In terms of Fourier series
\[
d\om = \sum_{n=-\infty}^{+\infty} d\om_n (t, x_2) e^{inx_1}, \quad d\om_{-n} = \overline{d\om_n},
\]
we have
\[
\pa_t d\om_n + inx_2 d\om_n = \frac{1}{Re} (\pa^2_{x_2} d\om_n - n^2 d\om_n ).
\]
Let
\[
d\om_n = d\Om_n e^{-inx_2t}, \quad d\Om_{-n} = \overline{d\Om_n},
\]
then
\[
\pa_t d\Om_n = \frac{1}{Re} [\pa^2_{x_2} d\Om_n - i2nt \pa_{x_2} d\Om_n - n^2 (t^2+1) d\Om_n ].
\]
When $Re = \infty$, i.e. for inviscid linearized Couette flow,
\[
d\Om_n (t, x_2) = d\Om_n (0, x_2).
\]
Thus \cite{Yud00}
\[
d\om (t, x_1, x_2) = d\om (0, x_1-x_2t, x_2).
\]
When $Re < \infty$, using Fourier transform
\[
d\Om_n = \int_{-\infty}^{+\infty} d\Om_{n\xi}(t) e^{i\xi x_2} d\xi , \quad d\Om_{(-n)\xi}
= \overline{d\Om_{n\xi}},
\]
we have
\[
\pa_t d\Om_{n\xi} = \frac{1}{Re} [-\xi^2 +2nt\xi - n^2 (t^2+1)] d\Om_{n\xi}.
\]
Thus
\[
d\Om_{n\xi}(t) = d\Om_{n\xi}(0) e^{-\frac{t}{Re} [(\xi - \frac{1}{2} nt)^2 + n^2
(\frac{1}{12}t^2+1)]},
\]
\[
d\Om_n (t, x_2) = \int_{-\infty}^{+\infty} d\Om_{n\xi}(0) e^{-\frac{t}{Re}
[(\xi - \frac{1}{2} nt)^2 + n^2 (\frac{1}{12}t^2+1)]} e^{i\xi x_2} d\xi ,
\]
and
\begin{eqnarray}
d\om &=& \sum_{n=-\infty}^{+\infty} e^{inx_1} e^{-inx_2t} d\Om_n (t, x_2) \non \\
&=& \sum_{n=-\infty}^{+\infty} e^{inx_1} e^{-inx_2t} \int_{-\infty}^{+\infty}
d\Om_{n\xi}(0) e^{-\frac{t}{Re}
[(\xi - \frac{1}{2} nt)^2 + n^2 (\frac{1}{12}t^2+1)]} e^{i\xi x_2} d\xi . \label{CouS}
\end{eqnarray}
Thus
\[
\| d\om \|_{H^0}^2 = 4\pi^2 \sum_{n=-\infty}^{+\infty} \int_{-\infty}^{+\infty}
|d\Om_{n\xi}(0)|^2 e^{-\frac{2t}{Re}
[(\xi - \frac{1}{2} nt)^2 + n^2 (\frac{1}{12}t^2+1)]} d\xi .
\]
When $Re = \infty$,
\[
\| d\om (t)\|_{H^0}^2 = \| d\om (0)\|_{H^0}^2 .
\]
When $Re < \infty$,
\[
\| d\om (t)\|_{H^0}^2 \leq \| d\om (0)\|_{H^0}^2 .
\]
Based on calculations such as
\[
\pa_{x_2} d\om = \sum_{n=-\infty}^{+\infty}[(-int)d\Om_n + \pa_{x_2} d\Om_n ] 
e^{inx_1} e^{-inx_2t} ,
\]
one can see that when  $Re = \infty$, $\| d\om (t)\|_{H^k}$ is finite but grows in time
as $t^k$ for $k = 1,2, \cdots $. As in the case of Example 1, the norm of the formal
derivative here is bounded, but the derivative does not exist due to the violation of
the definition of Fr\'echet derivative when $Re = \infty$. When $Re \ra \infty$, the norm of the
derivative in the Navier-Stokes case approaches the norm of the formal
derivative in the Euler case. Since the linear Couette shear is in $C^\infty$, it is not a base solution 
that is in $H^k$  but not in $H^{k+1}$. Historically linear Couette shear
was regarded as one the canonical flows for the study of transition to turbulence. Now we
understand that it is not a good representative of transition to turbulence. Indeed, it is both
linearly and nonlinearly stable for all values of Reynolds number including infinity \cite{Rom73}.
On the other hand, transition to turbulence does happen in lab Couette flow. Our explanation
of such a transtion is that states arbitrarily close to the linear Couette shear are linearly
unstable \cite{LL11}. In a typical transition to turbulence, the rough dependence may play a
significant role when the Reynolds number is large enough \cite{Li18}. But for linear Couette
flow, we believe that the transition is due to what we just mentioned \cite{LL11}. The base solution of our next example
does satisfy the criterion of being in $H^k$  but not in $H^{k+1}$. 

\section{Example 3 - periodic boundary condition}

For 2D Navier-Stokes equations under periodic boundary condition with period domain
$[0, 2\pi]\times [0, 2\pi]$, we have the one-parameter family of exact solutions \cite{Li17},
\[
u_1 = \sum_{n=1}^{\infty} \frac{1}{n^{3+\ga}}e^{-\frac{n^2 t}{Re}} \sin [ n (x_2 - \sg t)],
\quad  u_2 = \sg ,
\]
which is a solution in the space $C^0([0, \infty ), H^3)$ for all values of the Reynolds number
including infinity, $\frac{1}{2} < \ga \leq 1$, and $\sg$ is the real parameter. The directional derivative in $\sg$ of the solution
operator along the above exact solutions $\pa_\sg F^t$ is given by \cite{Li17}
\[
\pa_\sg u_1 = \sum_{n=1}^{\infty} \frac{-t}{n^{2+\ga}}e^{-\frac{n^2 t}{Re}} \cos [ n (x_2 - \sg t)],
\quad  \pa_\sg u_2 = 1.
\]
The norm of the derivative of the solution operator $\pa_\sg F^t$ has the lower bound \cite{Li17},
\[
\| \pa_\sg F^t \|_{H^3} > \sqrt{2} \pi  + \frac{\pi}{\sqrt{e}} t^\ga \left (\frac{\sqrt{t}\sqrt{Re}}
{2\sqrt{2}}\right )^{1-\ga} .
\]
As $Re \ra \infty$,
\[
\| \pa_\sg F^t \|_{H^3} \ra \infty ,
\]
thus
\[
\| \na F^t \|_{H^3} \ra \infty ,
\]
since $\| \na F^t \|_{H^3} \geq \| \pa_\sg F^t \|_{H^3}$. When $Re = \infty$, direct calculation
shows that indeed
\[
\| \pa_\sg F^t \|_{H^3} = \infty ,
\]
thus
\[
\| \na F^t \|_{H^3} = \infty .
\]
This example shows that the norm of the derivative of the solution
operator along the family of exact solutions is infinite in the Euler case, and approaches
infinity in the Navier-Stokes case as the Reynolds number approaches infinity.

\section{The generic case - numerical simulations}

\begin{figure}[ht] 
\centering
\subfigure[$Re=1000$]{\includegraphics[width=2.3in,height=2.3in]{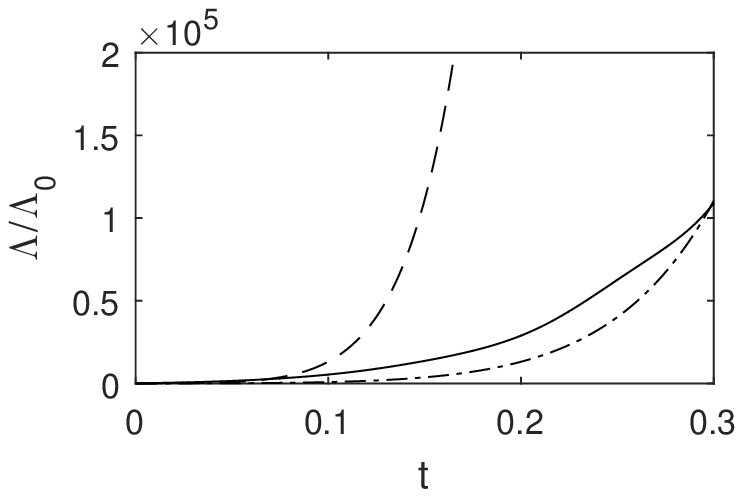}}
\subfigure[$Re=1000$]{\includegraphics[width=2.3in,height=2.3in]{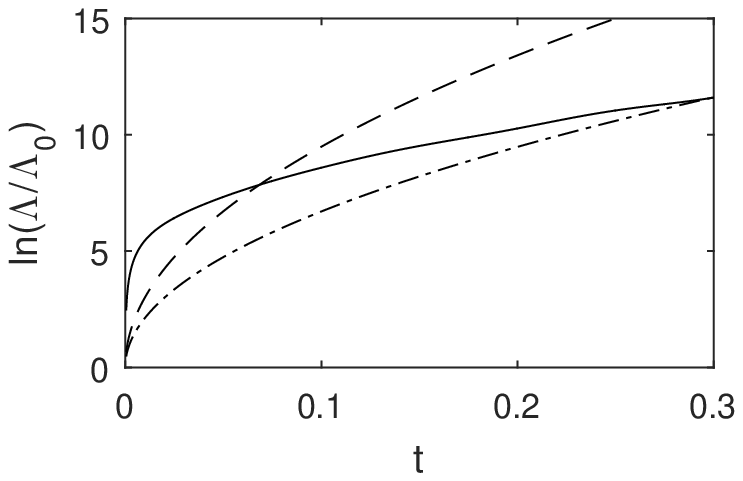}}
\caption{The solid curve is the numerical result of the super fast growth of perturbations where $\La (t) = \| du(t) \|_{H^3}$. 
The lower fitting dashed curve is $e^{21.2 \sqrt{t}}$. The closest fitting dashed curve is $e^{30 \sqrt{t}}$.}
\label{NF}
\end{figure}

We conducted extensive numerical simulations under periodic boundary condition \cite{FL19}\cite{LHBF19}. At high
Reynolds number (even moderate Reynolds number), along generic solutions, generic perturbations
amplify superfast (i.e. faster than exponential growth), typically as shown in Figure \ref{NF}.
In Figure \ref{NF}(a), we plot the growth of the $H^3$ norm of the perturbation under the linearized
Navier-Stokes dynamics at Reynolds number $Re=1000$, and in Figure \ref{NF}(b), we plot the same
figure in vertical $\ln$-scale. Clearly, the amplification of the perturbation is faster than
exponential growth. We want to emphasize that such superfast amplification is observed even for base solutions in $H^4$ (in fact $C^\infty$), and for
moderate Reynolds number.

\begin{figure}[ht]
\centering
\includegraphics[width=3.25in,height=2.25in]{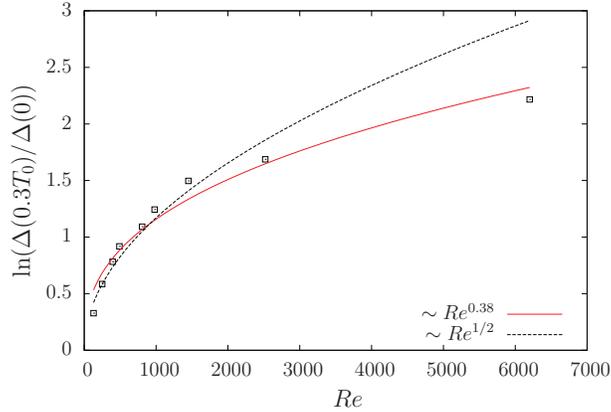}
\caption{The Reynolds number dependence of the amplifications of perturbations
for different Reynolds up to the time $t = 0.3 T_0$, where $\Dl$ is the $H^0$ norm of the perturbation, and  $T_0$ is the large 
eddy turnover time
which is around $2$. The dashed curve is a fit to $\sim \sqrt{Re}$, and the red (grey) curve
is a fit to $\sim Re^{0.38}$.}
\label{ARe}
\end{figure}

In \cite{LHBF19}, we conduct a large direct numerical simulation on the 3D Navier-Stokes equations
with resolution up to $2048^3$ and Reynolds number up to $6210$. We first run the numerical simulation for a long time 
until the dynamics reaches homogeneous and isotropic turbulence. The vorticity in this homogeneous and isotropic turbulence
is very large, and we use $H^0$ norm to measure the perturbations. We are simulating the situation that the $H^1$ norm 
of the base solutions are very large, and we measure their $H^0$ norm. The square root nature of both time and Reynolds 
number in the exponent of (\ref{UB}) is verified. Figure \ref{ARe} 
shows the Reynolds number dependence of the amplifications of perturbations. 
Thus we numerically verified our theory that fully developed turbulence is
initiated, developed and maintained by such superfast growth of ever existing perturbations!


\begin{thebibliography}{99}

\bibitem{FL19} Z. Feng, Y. Li, Short term unpredictability of high Reynolds number turbulence - rough dependence on initial data, {\it Submitted} (2019).

\bibitem{Inc15} H. Inci, On the regularity of the solution map of the incompressible Euler equation, {\it Dyn. Partial Differ. Equ.} {\bf 12, no.2} (2015), 97-113.

\bibitem{IL18} H. Inci, Y. Li, Nowhere-differentiability of the solution map of 2D Euler equations on bounded spatial domain, {\it Dyn. Partial Differ. Equ.} {\bf 16, no.4} (2019), 383-392.

\bibitem{Kat72} T. Kato, Nonstationary flows of viscous and ideal fluids in $R^3$, {\it J. Funct. 
Anal.} {\bf 9} (1972), 296-305.

\bibitem{Kat75} T. Kato, Quasi-linear equations of evolution, with applications to partial differential 
equations, {\it Lect. Notes in Math., Springer} {\bf 448} (1975), 25-70.

\bibitem{Li14} Y. Li, The distinction of turbulence from chaos --- rough dependence on initial data, {\it Electronic Journal of Differential Equations} {\bf 2014, no.104} (2014), 1-8.

\bibitem{Li17} Y. Li, Rough dependence upon initial data exemplified by explicit solutions and the effect of viscosity, {\it Nonlinearity} {\bf 30} (2017), 1097-1108.

\bibitem{Li18} Y. Li, Linear hydrodynamic stability, {\it Notices of the AMS} {\bf 65, no.10} (2018), 1255-1259.

\bibitem{LHBF19} Y. Li, R. Ho, A. Berera, Z. Feng, Superfast amplication and superfast nonlinear saturation of perturbations as the
mechanism of turbulence, {\it arXiv:1908.04838}
 
\bibitem{LL11} Y. Li, Z. Lin, A resolution of the Sommerfeld paradox, {\it SIAM J. Math. Anal.} {\bf 43, no.4} (2011), 1923-1954.

\bibitem{Rom73} V. Romanov, Stability of plane-parallel Couette flow, {\it Funct. Anal. Appl.}
{\bf 7} (1973), 137–146.

\bibitem{Yud00} V. Yudovich, On the loss of smoothness of the solutions of the Euler equations and the
inherent instability of flows of an ideal fluid, {\it Chaos} {\bf 10, no.3} (2000), 705-719.


\end{thebibliography}
\end{document}